\documentclass{article}
\usepackage{amsmath,amssymb,amsfonts,latexsym,stmaryrd,diagrams,mathrsfs}

\newtheorem{thm}{Theorem}[section]
\newtheorem{lem}{Lemma}[section]

\newtheorem{defn}{Definition}[section]

\newenvironment{prf}{\noindent{\em Proof: }}
               {\nopagebreak\hfill$\Box$\vspace{2ex}}

\newcommand{\arr}{\rightarrow}
\newcommand{\Pow}{\mathcal{P}}

\newcommand{\Sub}{{\mathsf{Sub}}}

\newcommand{\N}{\mathbb{N}}

\def\lbr{\mathopen{{[\kern-0.14em[}}}   
\def\rbr{\mathclose{{]\kern-0.14em]}}}  

\newcommand{\ev}{\mathsf{ev}}

\newcommand{\inc}{\hookrightarrow}

\newcommand{\strarr}{\mbox{$\circ \kern-0.4em \arr$}}
\newcommand{\eps}{\,\varepsilon\,}
\newcommand{\neps}{\; \slash \kern-0.7em \eps \,}

\newcommand{\RT}{{\mathbf{RT}}}

\newcommand{\A}{{\mathcal{A}}}

\newcommand{\C}{{\mathcal{C}}}

\newcommand{\E}{{\mathcal{E}}}
\newcommand{\F}{{\mathcal{F}}}

\newcommand{\opp}{{\mathsf{op}}}

\newcommand{\forget}[1]{}

\newcommand{\mono}{\rightarrowtail}

\newcommand{\cov}{\rightarrowtriangle}

\newcommand{\sheafify}{{\mathsf{a}}}
\newcommand{\incl}{{\mathsf{i}}}
\newcommand{\Se}{\mathbf{Set}}

\newcommand{\IZF}{\mathsf{IZF}}

\newcommand{\Sm}{{\mathcal{S}}}

\newcommand{\Eff}{\mathcal{E}\kern-0.14em\mathit{ff}}

\newcommand{\Asm}{{\mathbf{Asm}}}
\newcommand{\Mod}{{\mathbf{Mod}}}

\newcommand{\Kcal}{{\mathcal{K}}}

\newcommand{\card}{{\mathsf{card}}}

\newcommand{\schloop}{{< \kern-0.40em  <}}

\newcommand{\eq}{{\mathit{eq}}}

\newcommand{\Trip}{{\mathscr{P}}}
\newcommand{\Prop}{{\mathit{Prop}}}
\newcommand{\Truth}{{\mathit{Tr}}}

\diagramstyle[PostScript=dvips]

\newarrow{Dashto}{}{dash}{}{dash}>
\newarrow{Embed}>--->
\newarrow{Into} C--->
\newarrow{To} ---->
\newarrow{Line} -----
\newarrow{Onto} ----{>>}
\newarrow{Cover} ----{triangle}
\newarrow{Equal} =====
\newarrow{Distto}--+->
\newarrow{Lin}----o

\begin{document}

\title{Models of Intuitionistic Set Theory \\ 
       in Subtoposes of Nested Realizability Toposes}

\author{S.~Maschio \\
Dipartimento di Matematica, Universit\`a di Padova \\
Via Trieste, Padova\\
{\small\tt samuele.maschio@math.unipd.it} \\[5mm]
T.~Streicher \\
Fachbereich 4 Mathematik, TU Darmstadt \\
Schlo{\ss}gartenstr.~7, D-64289 Darmstadt, Germany\\
{\small\tt streicher@mathematik.tu-darmstadt.de}}
\date{\today}

\maketitle
\section*{Introduction}

Given a \emph{partial combinatory algebra} (pca) $\A$ (see e.g.\ 
\cite{VOO08}) together with a subpca $\A_\#$ of $\A$ we will construct
the \emph{nested realizability topos} $\RT(\A,\A_\#)$ as described 
in \cite{BVO02} (without giving it a proper name there). It is well known 
(from e.g.\ \cite{VOO08}) that $\RT(\A,\A_\#)$ appears as the exact/regular 
completion of its subcategory $\Asm(\A,\A_\#)$ of assemblies. 
In \cite{BVO02} the authors considered two complementary subtoposes of 
$\RT(\A,\A_\#)$, namely the \emph{relative realizability topos} 
$\RT_r(\A,\A_\#)$ and the \emph{modified relative realizability topos}
$\RT_m(\A,\A_\#)$, respectively. 

Within nested realizability toposes we will identify a class of 
\emph{small maps} giving rise to a model of intuitionistic set theory $\IZF$ 
(see \cite{FRI73,MCC84}) as described in \cite{JM95}. For this purpose we
first identify a class of display maps in $\Asm(\A,\A_\#)$ which using a 
result of \cite{VDBM08} gives rise to the desired class of small maps in
the exact/regular completion $\RT(\A,\A_\#)$ of $\Asm(\A,\A_\#)$.

For showing that the subtoposes $\RT_r(\A,\A_\#)$ and  $\RT_m(\A,\A_\#)$
also give rise to models of $\IZF$ we will prove the following general result.
If $\E$ is a topos with a class $\Sm$ of small maps and $\F$ is a subtopos
of $\E$ then there is a class $\Sm_\F$ of small maps in $\F$ which is obtained
by closing sheafifications of maps in $\Sm$ under quotients in $\F$.

As explained in subsections \ref{modreal} and \ref{HerbTop} below this
covers also the Modified Realizability topos as studied in \cite{VOO97} and 
the more recent Herbrand topos of van den Berg.

\section{Nested Realizability Toposes and some of their Subtoposes}

Given a pca $\A$ in an elementary topos $\mathscr{S}$ we may construct the 
realizability topos $\RT_{\mathscr{S}}(\A)$ relative to $\Sm$ as described in
\cite{VOO08}. If $\mathscr{S}$ is the Sierpi\'nski topos $\Se^{\mathbf{2}^\opp}$
then a ``nested pca'', i.e.\ a pca $\A$ together with a subpca $\A_\#$ gives
rise to a pca internal to  $\Se^{\mathbf{2}^\opp}$ from which one may construct
the ``nested realizability topos'' $\RT(\A,\A_\#)$ as described in 
\cite{BVO02,VOO08}.\footnote{In \cite{BVO02} they do not give a name
to this topos and, moreover, write $\RT(\A,\A_\#)$ for the relative 
realizability subtopos of the nested realizability topos.}
 Within $\RT(\A,\A_\#)$ there is a unique nontrivial subterminal object $u$
giving rise to the \emph{open} subtopos induced by the closure operator 
$u \to (-)$ and the complementary subtopos induced by the closure
operator $u \vee (-)$ as described in \cite{BVO02}.

Next we will give more elementary descriptions of $\RT(\A,\A_\#)$
and the above mentioned subtoposes.

\subsection{The Nested Realizability Topos $\RT(\A,\A_\#)$}

Let $\A$ be a pca whose partial application is denoted by juxtaposition
and $\A_\#$ be a subpca of $\A$, i.e.\ $\A_\#$ is a subset of $\A$ closed 
under application and there are elements $k$ and $s$ of $\A_\sharp$ such that
for all $x,y,z \in \A$ it holds that $kxy = x$, $sxyz \simeq xz(yz)$ and
$sxy$ is always defined. We write $i$ for $skk$ and $\bar{k}$ for $ki$ which,
obviously, satisfy the equations $ix =  x$ and $\bar{k}xy = y$, respectively.
We write $p$, $p_0$ and $p_1$ for elements of $\A$ such that $px_0x_1$ is 
always defined and $p_i(px_0x_1) = x_i$ for $i = 0,1$.
For every natural number $n$ we write $\underline{n}$ for the corresponding 
numeral as defined in \cite{VOO08}. Notice that $k$, $\bar{k}$, $p$, $p_0$, 
$p_1$ and the numerals $\underline{n}$ are all elements of $\A_\#$. 

Since subsets of $\A$ are the propositions of the realizability topos $\RT(\A)$
it is useful to fix some notation for the propositional connectives 
\begin{enumerate}
\item[] $A \rightarrow B = 
         \left\{ a \in \A \mid  ax \in B \textrm{ for all } x \in A \right\}$
\item[] $A \wedge B = \left\{ pxy \mid x\in A,  y \in B \right\}$
\item[] $A \vee B = \left(\{ k \} \wedge A\right) \cup  
                    \left(\{ \bar{k} \} \wedge B \right)$
\end{enumerate}
Propositions of the nested realizability topos $\RT(\A,\A_\#)$ will be pairs
$A = (A_p,A_a) \in \Pow(\A) \times \Pow(\A_\#)$ such that $A_a \subseteq A_p$
where we call $A_p$ and $A_a$ the set of \emph{potential} and \emph{actual}
realizers, respectively. We write $\Sigma(\A,\A_\#)$ for the set of these
propositions. The above notation for propositional connectives
is adapted to the current class of propositions as follows
\begin{enumerate}
\item[]  $A \rightarrow B = \left( A_p \rightarrow B_p, \A_\# \cap (A_p \rightarrow B_p) \cap (A_a \rightarrow B_a)\right)$
\item[] $A \wedge B = (A_p \wedge B_p, A_a \wedge B_a)$
\item[] $A \vee B = (A_p \vee B_p, A_a \vee B_a)$
\end{enumerate}

For the realizability tripos $\mathscr{P}(\A)$ induced by the pca $\A$ see
\cite{VOO08}. The nested realizability tripos $\Trip(\A,\A_\#)$ over $\Se$ 
induced by the nested pca $\A_\# \subseteq \A$ is defined as follows. 
For a set $I$ the fibre $\Trip(\A,\A_\#)(I)$ is given by the set 
$\Sigma(\A,\A_\#)^I$ preordered by the relation $\vdash_I$ defined as
\[ \phi \vdash_{I} \psi \qquad \textrm{ if and only if } \qquad 
   \bigcap_{i\in I}(\phi(i) \rightarrow \psi(i))_{a} \neq \emptyset \]
for $\phi,\psi \in \Trip(\A,\A_\#)(I)$.
For $u : J \to I$ reindexing along $u$ is given by precomposition with $u$
and denoted as $u^*$.
The fibres are preHeyting algebras where the propositional connectives are
given by applying the operations $\rightarrow$, $\wedge$ and $\vee$ pointwise. 
It is easy to check that $u^*$ commutes with the propositional connectives 
in the fibres. For a map $u : J \to I$, the reindexing $u^*$ has left and 
right adjoints $\exists_u$ and $\forall_u$, respectively, given by
\[ \exists_u(\phi)(i) = 
   \left(\bigcup_{u(j)=i} \phi_p(j), \bigcup_{u(j)=i}\phi_a(j)\right) \]
\[ \forall_u(\phi)(i) = \left(\bigcap_{j \in J}(Eq(u(j),i) \rightarrow 
   \phi(j))_p, \bigcap_{j \in J}(Eq(u(j),i) \rightarrow 
   \phi(j))_a \right) \]
where $Eq(x,y) = \left(\left\{a \in \A \mid x = y\right\},\left\{a \in \A_\# \mid x = y\right\}\right)$. It is straightforward to check that the so defined
quantifiers satisfy the respective Beck-Chevalley conditions. 
The identity on $\Sigma(\A,\A_\#)$ gives rise to a generic family and,
therefore, the fibered preorder $\Trip(\A,\A_\#)$ is actually a tripos
in the sense of \cite{HJP80}. 

We write $\RT(\A,\A_\#)$ for the ensuing topos.
 
\subsection{Some Subtoposes of $\RT(\A,\A_\#)$}\label{subtopex}
 
In $\RT(\A,\A_\#)$ there is a nontrivial subterminal $u = (\A,\emptyset)$
giving rise to two complementary subtoposes induced by the closure operators
$o_u(p) = u \to p$ and $c_u(p) = u \vee p$ as in \cite{BVO02}. We denote the 
open subtopos induced by $o_u$ by $\RT_r(\A,\A_\#)$ and the complementary 
subtopos induced by $c_u$ by $\RT_m(\A,\A_\#)$. In \cite{BVO02} these two
subtoposes are referred to as the \emph{relative} and the \emph{modified 
relative realizability} topos, respectively.

For sake of concreteness and later reference in the following two subsections
we give an elementary and explicit construction of triposes inducing 
$\RT_r(\A,\A_\#)$ and $\RT_m(\A,\A_\#)$, respectively.

\subsubsection{The Relative Realizability Topos $\RT_r(\A,\A_\#)$}

is induced by the tripos $\Trip_r(\A,\A_\#)$ over $\Se$ which we describe next.
Let $\Sigma_r(\A,\A_\#) = \Pow(\A)$. The fibre of $\Trip_r(\A,\A_\#)$ over $I$
is given by the preorder $\left(\Pow(\A)^{I},\vdash^r_I\right)$ where
\[ \phi \vdash_I^r \psi \qquad \textrm{ if and only if } \qquad
   \A_\# \cap \bigcap_{i\in I} (\phi(i) \rightarrow \psi(i)) \neq \emptyset \] 
and as usual reindexing is given by precomposition. At first sight
this tripos looks like the tripos $\Trip(\A)$ inducing the realizability topos
$\RT(\A)$ but notice that entailment in the fibres is defined in a more 
restrictive way, namely by requiring that the entailment be realized by an 
element of $\A_\#$ and not just an element of $\A$. Nevertheless, the
propositional connectives, quantifiers and the generic family of 
$\RT_r(\A,\A_\#)$ can be constructed according to the same recipes as for
$\Trip(\A)$ (see \cite{VOO08}). 

There is an obvious logical morphism from $\Trip_r(\A,\A_\#)$ to $\Trip(\A)$
which is the identity on objects. But there is also an injective geometric
morphism from $\Trip_r(\A,\A_\#)$ to $\Trip(\A,\A_\#)$ sending a family 
$\phi \in \Pow(\A)^I$ to the family 
$\lambda i{:}I. (\phi(i),\A_\# \cap \phi(i))$. 
These morphisms between triposes over $\Se$ extend to morphisms between the 
associated toposes as described in \cite{VOO08}.

\subsubsection{The Modified Relative Realizability Topos $\RT_m(\A,\A_\#)$}\label{modreal}

is induced by the tripos $\Trip_m(\A,\A_\#)$ over $\Se$ which is obtained
from $\Trip(\A,\A_\#)$ by restricting the fibre $\Trip(\A,\A_\#)(I)$ to the 
set of all $\phi \in \Sigma(\A,\A_\#)^I$ with
$\A_\# \cap \bigcap_{i \in I} \phi_p(i) \neq \emptyset$.
The logical structure is essentially inherited from $\Trip(\A,\A_\#)$ though
now and then one has to insert the closure operator $c_U$ in order to stay
within $\Trip_m(\A,\A_\#)$.
A generic family for $\Trip_m(\A,\A_\#)$ is given by the identity on
$\Sigma_m(\A,\A_\#) = \{ u \vee p \mid p \in \Sigma(\A,\A_\#) \}$.
The obvious inclusion of $\Trip_m(\A,\A_\#)$ into  $\Trip(\A,\A_\#)$ gives
rise to the inclusion of $\RT_m(\A,\A_\#)$ into $\RT(\A,\A_\#)$.

Notice that in case $\A = \A_\#$ we obtain the \emph{modified realizability}
topos as described in \cite{VOO97,VOO08} for the case where $\A$ is the
first Kleene algebra $\mathcal{K}_1$.

Another prominent example of a modified relative realizability model can
be found in a paper by J.~R.~Moschovakis \cite{MOS} from 1971 where she
constructed a model for a theory $\mathrm{INT}$ of Brouwerian intuitionism 
validating the proposition that all functions on natural numbers are not not 
recursive, i.e.\
that there are no non-recursive functions on the natural numbers. Of course,
the model of \cite{MOS} was not constructed in topos-theoretic terms but it
is equivalent to the interpretation of the system considered in 
\emph{loc.cit.}\ in the topos $\RT_m(\mathcal{K}_2,\mathcal{K}_2^{\mathit{rec}})$
where $\mathcal{K}_2$ is the second Kleene algebra whose underlying set is 
Baire space $\N^\N$ and $\mathcal{K}_2^{\mathit{rec}})$ is the sub-pca of 
recursive sequences of natural numbers. The ensuing interpretation of 
$\mathrm{INT}$ was called \emph{$G$-realizability} in \emph{loc.cit.}

\subsubsection{The Herbrand Realizability Topos}\label{HerbTop}

As shown by J.~van Oosten, see Lemma~3.2 of \cite{JohGl}, B.~van den Berg's 
\emph{Herbrand realizability topos} over a pca $\A$ arises as a subtopos 
of $\RT(\A,\A)$ induced by some closure operator on $\Trip(\A,\A)$.
Moreover, as shown in \emph{loc.cit.}\ it is disjoint from the open subtopos 
$\RT_r(\A,\A)$ equivalent to $\RT(\A)$.

\subsection{Assemblies induced by $\Trip(\A,\A_\#)$}

As described in \cite{VOO08} for every tripos $\Trip$ (over $\Se$) one may 
consider the full subcategory $\Asm(\Trip)$ of \emph{assemblies} in
$\Se(\Trip)$, i.e.\ subobjects of objects of the form $\Delta(S)$ where
$S \in \Se$ and $\Delta : \Se \to \Se(\Trip)$ is the \emph{constant objects}
functor sending a set $S$ to $(S,\exists_{\delta_S}(\top_S))$.\footnote{In \cite{VOO08} the constant objects functor is denoted by $\nabla$ because in
case of realizability triposes it is right adjoint to the global elements 
functor $\Gamma$. However, in case of triposes induced by a complete Heyting 
algebra the constant objects functor is left adjoint to $\Gamma$. However, there 
are also triposes where the constant objects functor is neither left nor right 
adjoint to $\Gamma$. We prefer the notation $\Delta$ since 
$\eq_S = \exists_{\delta_S}(\top_S)$ is the (Lawvere) equality predicate 
on the set $S$ in the sense of the tripos $\Trip$.}

One can show that the category $\Asm(\Trip(\A,\A_\#))$ is equivalent to
the category $\Asm(\A,\A_\#)$ whose objects are pairs $X = (|X|,E_X)$
where $|X|$ is a set and $E_X : |X| \to \Sigma(\A,\A_\#)$ with $E_X(x)_p
\neq \emptyset$ for all $x \in |X|$. An arrow from $X$ to $Y$ is a function 
$f : |X| \rightarrow |Y|$ such that $E_X \vdash_{|X|} f^*E_Y$.

As follows from \cite{VOO08} Cor.~2.4.5 the topos $\RT(\A,\A_\#)$ appears
as the exact/regular completion of $\Asm(\A,\A_\#)$.

For further reference we note the following

\begin{thm}\label{AsmplcH}
$\Asm(\A,\A_\#)$ is a locally cartesian closed Heyting category
with stable and disjoint finite sums with a generic monomorphism
$\top : \Truth \mono \Prop$.\footnote{``generic'' means that all monos 
can be obtained as pullbacks of $\top : \Truth \mono \Prop$
but we may have $f^*\top \cong g^*\top$ for different $f$ and $g$}
\end{thm}
\begin{prf}
The locally cartesian closed structure is constructed as in the case of
$\Asm(\A)$, i.e.\ assemblies within $\RT(\A)$ where $\A$ is a pca. Similarly,
one shows that $\Asm(\A,\A_\#)$ is a Heyting category and it has stable and
disjoint finite sums.

Finally we exhibit a generic mono $\top : \Truth \mono \Prop$. The object
$\Prop$ is defined as $\Delta(\Sigma(\A,\A_\#))$. The underlying set of $\Truth$
is the subset of $\Sigma(\A,\A_\#)$ consisting of those pairs $A = (A_p,A_a)$
where $A_p \neq \emptyset$ and $E_\Truth(A) = A$.
\end{prf}

Notice, however, that in general $\Asm(\A,\A_\#)$ is not well-pointed.

\section{Some Facts about Small Maps}

A Heyting category is a regular category $\C$ where for all $f : Y \to X$ 
in $\C$ the pullback functor $f^{-1} : \Sub_\C(X) \to \Sub_\C(Y)$ has a right 
adjoint $\forall_f$. It is a Heyting pretopos iff, moreover, it has stable 
disjoint finite sums and every equivalence relation is effective (i.e.\ 
appears as kernel pair of its coequalizer).

\begin{defn}\label{Sdef}
Let $\C$ be a locally cartesian Heyting category 
with stable and disjoint finite sums and a natural numbers object $N$.
For a class $\Sm$ of maps in $\C$ we consider the following properties.
\begin{itemize}
\item[{\bf (A0)}] (Pullback Stability) For a pullback square
\begin{diagram}[small]
D \SEpbk & \rTo^q & B \\
\dTo^g & & \dTo_f \\
C & \rTo_p & A
\end{diagram}
in $\C$ from $f \in \Sm$ it follows that $g \in \Sm$.
\item[{\bf (A1)}] (Descent) If in a pullback square as above $p$ is a cover,
i.e.\ a regular epimorphism, then $f \in \Sm$ whenever $g \in \Sm$.
\item[{\bf (A2)}] (Sums) If $f$ and $g$ are in $\Sm$ then $f + g$ is in $\Sm$.
\item[{\bf (A3)}] (Finiteness) The maps $0 \to 1$, $1 \to 1$ and $1 + 1 \to 1$ 
are in $\Sm$.
\item[{\bf (A4)}] (Composition) Maps in $\Sm$ are closed under composition.
\item[{\bf (A5)}] (Quotient) If $f \circ e$ is in $\Sm$ and $e$ is a cover
                             then $f$ is in $\Sm$.
\item[{\bf (A6)}] (Collection) Any arrows $p : Y \to X$ and $f : X \to A$
where $p$ is a cover and $f \in \Sm$ fit into a quasipullback\footnote{A square
is a quasipullback if the mediating arrow to the pullback square is a cover.}
\begin{diagram}[small]
Z & \rTo & Y & \rCover^p & X \\
\dTo^g & & & &  \dTo_f \\
B & & \rCover_h & & A
\end{diagram}
where $g \in \Sm$ and $h$ is a cover.
\item[{\bf (A7)}] (Representability) There is a universal family
$\pi : E \to U$ in $\Sm$ 
such that every $f : Y \to X$ in $\Sm$ fits into a diagram
\begin{diagram}
Y & \lCover & Y^\prime & \SEpbk \rTo & E \\
\dTo^f & \mathrm{qpb} & \dTo^{f^\prime} & & \dTo_\pi \\
X & \lCover & X^\prime & \rTo & U \\
\end{diagram} 
where the left square is a quasipullback and the right square is a pullback.
\item[{\bf (A8)}] (Infinity) The terminal projection $N \to 1$ is in $\Sm$.
\item[{\bf (A9)}] (Separation) All monomorphisms are in $\Sm$.
\end{itemize}
A class $\Sm$ of maps in $\C$ validating properties {\bf (A0)}--{\bf (A9)}
is called a class of small maps.
\end{defn}

The following theorem will be essential later on.

\begin{thm}\label{bennoextended}
Let $\C$ be a Heyting category with stable and disjoint finite sums and $\Sm$
be a class of small maps in $\C$. Let $\bar{\C}$ be the exact/regular completion
of $\C$ and $\bar{\Sm}$ the class of maps $f$ in $\bar{\C}$ which fit into a
quasipullback
\begin{diagram}[small]
\cdot & \rCover & \cdot \\
\dTo^g & & \dTo_f \\
\cdot & \rCover & \cdot 
\end{diagram} 
with $g$ in the subcategory $\C$ of $\bar{\C}$.

Then $\bar{\Sm}$ is a class of small maps 
within the Heyting pretopos $\bar{\C}$.
\end{thm}
\begin{prf}
That $\bar{\Sm}$ validates conditions {\bf (A0)}--{\bf (A8)} follows from
Lemma~5.8 and Propositions 6.2 and 6.21 in \cite{VDBM08}. 

Condition {\bf (A9)} holds for $\bar{\Sm}$ in $\bar{\C}$ for the following 
reason. Let $m : B \mono A$ be a mono in 
$\bar{\C}$. Since $\bar{\C}$ is the exact completion of $\C$ there is a
cover $p : X \cov A$ with $X$ in $\C$. Then for the pullback
\begin{diagram}[small]
Y \SEpbk & \rCover^q & B \\
\dEmbed^n & & \dEmbed_m \\
X & \rCover_p & A
\end{diagram}
in $\bar{\C}$ we know that $q$ is a cover and $n$ is a mono. It follows
from Lemma~2.4.4 of \cite{VOO08} that $Y$ is isomorphic to an object in $\C$.
\end{prf}

 \section{Small Maps  in Nested Realizability Toposes}

We will first identify within $\Asm(\A,\A_\#)$ a class $\Sm$ of small maps 
so that we can apply Theorem~\ref{bennoextended} to it in order to obtain a
class $\bar{\Sm}$ of small maps on $\RT(\A,\A_\#)$ which is known to arise 
as the exact/regular completion of $\Asm(\A,\A_\#)$ (see section 2.4 of 
\cite{VOO08} for more details). 

However, for showing that $\bar{\Sm}$ is closed under power types we have to
appeal to Lemma~27 of \cite{VDBM11} guaranteeing that if $\Asm(\A,\A_\#)$
has \emph{weak} power types under which $\Sm$ is closed then $\RT(\A,\A_\#)$
has power objects under which $\bar{\Sm}$ is closed.

\subsection{Small maps in $\Asm(\A,\A_\#)$}

For constructing a class of small maps in $\Asm(\A,\A_\#)$ let us first choose 
a \emph{strongly inaccessible} cardinal $\kappa$ exceeding the cardinality 
of $\A$.

\begin{thm}\label{Sclos1} 
Let $\Sm$ be the class of all maps $f:Y \to X$ in $\Asm(\A,\A_\#)$ 
such that $\card\left(f^{-1}(x)\right) < \kappa$ for all $x \in |X|$. 
Then $\Sm$ is a class of small maps in $\Asm(\A,\A_\#)$ 
in the sense of Def.~\ref{Sdef}.
\end{thm}
\begin{prf}
Conditions {\bf (A0)} and {\bf (A1)} follow from the fact that 
the forgetful functor from $\Asm(\A,\A_\#)$ to $\Se$ preserves finite limits 
and covers.

Since the forgetful functor from $\Asm(\A,\A_\#)$ to $\Se$ 
preserves finite sums condition {\bf (A2)} holds.

Since $\kappa$ is infinite all maps in $\Asm(\A,\A_\#)$ with finite fibres 
are in $\Sm$. For this reason {\bf (A3)} and {\bf (A9)} trivially hold.

Condition {\bf (A4)} holds since $\kappa$ is regular.

For {\bf (A5)} suppose $f \circ e$ is in $\Sm$ and $e$ is a cover. Then 
the fibres of $f$ have cardinalities $< \kappa$ since by assumption the 
fibres of $f \circ e$ have cardinalities $< \kappa$ and the underlying map 
of $e$ is onto.

Condition {\bf (A8)} holds since $\kappa$ exceeds the cardinality of $\N$.

For showing that {\bf (A6)} holds suppose $p : Y \to X$ is a cover and 
$f : X \to A$ is in $\Sm$. Since $p$ is a cover the underlying map of $p$ 
(also denoted by $p$) is onto and there exists $a \in \A_\#$ such that 
for all $x \in |X|$ it holds that
\begin{enumerate}
\item[(1p)] if $b \in E_X(x)_p$ then $ab{\downarrow}$ and $ab \in E_Y(y_{x,b})$
            for some $y_{x,b} \in p^{-1}(x)$ and
\item[(1a)] if $b \in E_X(x)_a$ then $ab \in E_Y(y_{x,b})_a$.
\end{enumerate}
Let $Z$ be the object of $\Asm(\A,\A_\#)$ whose underlying set 
$|Z| = \{ y_{x,b} \mid x \in |X|, b \in E_X(x)_p \}$ and 
$E_Z(y) = E_Y(y)$ for $y \in |Z|$. Let $i : Z \inc Y$
be the obvious inclusion of $Z$ into $Y$. Then the rectangle
\begin{diagram}
Z & \rEmbed^i & Y & \rCover^p & X \\
\dTo & & & &  \dTo_f \\
A & & \rEqual & & A
\end{diagram}
is a quasipullback since $p \circ i$ is a cover. Since the fibres of 
$p \circ i$ have cardinality $\leq \card(\A) < \kappa$ the map $p \circ i$ 
is in $\Sm$. Thus, by {\bf (A4)} the map $f \circ p \circ i : Z \to A$ is 
in $\Sm$, too.

Condition {\bf (A7)} holds in a very strong sense because we can exhibit 
a generic map $\pi : E \to U$ in $\Sm$, i.e.\ $\pi \in \Sm$ and all maps 
in $\Sm$ can be obtained as pullbacks of the generic map $\pi$. The codomain 
$U$ of $\pi$ is given by
\[ \Delta\left(\left\{ X \in \Asm(\A,\A_\#) \mid 
                   |X| \subseteq \kappa, \card(|X|) < \kappa\right\}\right) \]
and its domain $E$ has underlying set
\[ |E| = \left\{ (X,x) \mid X \in |U|, x \in |X| \right\} \]
and whose existence predicate is given by $E_E(X,x) = E_X(x)$. The map 
$\pi : E \to U$ is given by projection on the first component, i.e.\ $\pi(X,x) = X$. 
Obviously, the map $\pi$ has fibres of cardinality $< \kappa$ and we leave it 
as a straighforward exercise for the reader to show 
that every map in $\Sm$ can actually be obtained as pullback of $\pi$.
\end{prf}

It is easy to check that the class $\Sm$ in $\Asm(\A,\A_\#)$ is closed under
dependent products, i.e.\ $\Pi_f g \in \Sm$ whenever $f$ and $g$ are in $\Sm$.
As a consequence for $a : A \to I$ and $b : B \to I$ in $\Sm$ their exponential
in the fibre over $I$, i.e.\ $a \to_I b = \Pi_a a^* b$, is in $\Sm$, too.
Moreover, the generic mono $\top : \Truth \mono \Prop$ constructed in 
Theorem~\ref{AsmplcH} like all monos is also an element of $\Sm$. Moreover,
the terminal projection $\Prop \to 1$ is in $\Sm$, too, since the underlying 
set of $\Prop$ has cardinality $< \kappa$. Accordingly, the object $\Truth$
is small, too.

For every object $X$ in $\Asm(\A,\A_\#)$ we may construct a 
\emph{weak power object} $\ni^w_X \mono \Prop^X {\times} X$ as follows
\begin{diagram}[small]
\ni^w_X \SEpbk & \rTo & \Truth \\
\dEmbed & & \dEmbed_\top \\
\Prop^X {\times} X & \rTo_\ev & \Prop
\end{diagram}
where $\ev : \Prop^X \times X \to \Prop$ is the evaluation map. If $X$ is 
small, i.e.\ $X \to 1$ is in $\Sm$, i.e.\ $\card(X) < \kappa$, then $\Prop^X$
is small, too, since $\card\left(\Prop^X\right) \leq \card(\Prop)^{\card(X)} 
< \kappa$ because $\kappa$ is inaccessible and 
$\card(\Prop) , \card(X) < \kappa$.
Notice that this construction of weak power objects also works in all slices. 

For future reference we summarize these considerations in the following

\begin{thm}\label{weakpowerAsm}
The category $\Asm(\A,\A_\#)$ has weak power objects and $\Sm$ is closed
under weak power objects.
\end{thm}

\subsection{Small maps in $\RT(\A,\A_\#)$}

It is well known from \cite{VOO08} (section 2.4) that $\RT(\A,\A_\#)$ is the 
exact/regular completion of $\Asm(\A,\A_\#)$. Let $\bar{\Sm}$ be the class of
maps defined in Theorem~\ref{bennoextended}. Now we can show easily that

\begin{thm}\label{smallnesttop}
$\bar{\Sm}$ is a class of small maps in $\RT(\A,\A_\#)$ 
which is also closed under power objects and thus also under exponentiation.
\end{thm}
\begin{prf}
It is an immediate consequence of Theorem~\ref{bennoextended} and 
Theorem~\ref{Sclos1} that $\bar{\Sm}$ is a class of small maps in $\RT(\A,\A_\#)$.
From Lemma 27 of \cite{VDBM11} and our Theorem~\ref{weakpowerAsm} it follows
that $\bar{\Sm}$ is also closed under power objects. 
It is well known that closure under powerobjects and subobjects entails 
closure under exponentiation.
\end{prf}

As pointed out by J.~van Oosten in private communication there is a
logical functor $F : \RT(\A,\A_\sharp) \to \RT(\A)$ which just ``forgets
the actual realizers''. Already in \cite{JM95} there has been identified
for every strongly inaccessible cardinal a class of small maps in $\RT(\A)$
from which our class of small maps in $\RT(\A,\A_\#)$ can be obtained as
inverse image under $F$.

\subsection{A Model of $\IZF$ in $\RT(\A,\A_\#)$}

It follows from the previous Theorem~\ref{smallnesttop} and
Theorem~5.6 of \cite{JM95} that the class $\bar{\Sm}$ 
of small maps in $\RT(\A,\A_\#)$ gives rise to an 
``initial $\mathsf{ZF}$-algebra'' within $\RT(\A,\A_\#)$. 
Accordingly, the nested realizability topos 
$\RT(\A,\A_\#)$ hosts a model of $\IZF$.

It is an open question (raised by J.~van Oosten) whether the above mentioned 
logical functor $F : \RT(\A,\A_\#) \to \RT(\A)$ preserves the initial 
$\mathsf{ZF}$-algebras arising from the respective classes of small maps.

\section{Small Maps for Subtoposes of $\RT(\A,\A_\#)$}

In the previous section we have endowed the nested realizability topos 
$\RT(\A,\A_\#)$ with a class $\bar{\Sm}$ of small maps in such a way that
it gives rise to a model of $\IZF$ in the sense of Algebraic Set Theory 
as described in \cite{JM95}. In this section we show how to extend this 
result to subtoposes of $\RT(\A,\A_\#)$.

\subsection{Transferring Classes of Small Maps to Subtoposes}

Let $\E$ be an elementary topos and $\Sm$ a class of small maps in $\E$.
Let $\sheafify \dashv \incl : \F \inc \E$ be a subtopos of $\E$. 
W.l.o.g.\ we assume that $\F$ is closed under isomorphisms in $\E$ and 
that $\sheafify f = f$ for $f \in \F$. We want to endow $\F$ with a class 
$\Sm_\F$ of small maps such that $\sheafify : \E \to \F$ sends $\Sm$ to
$\Sm_\F$. Thus, it is tempting to define $\Sm_\F$ as $\sheafify\F$ by which
we denote the closure under isomorphism in $\F$ of the image of $\sheafify$. 
But then there are problems with condition {\bf (A5)} because epimorphisms 
in $\F$ need not be epic in $\E$. 
In order to overcome this problem we define $\Sm_\F$ as follows

\begin{defn}\label{SFdf} Let $\Sm_\F$ be the class of all maps $f : B \to A$
in $\F$ fitting into a quasipullback
\begin{diagram}[small]
\sheafify Y & \rTo &  B \\
\dTo^{\sheafify g} & \mathrm{qpb} & \dTo_f \\
\sheafify X & \rCover_e & A
\end{diagram}
in $\F$ for some $g : Y \to X$ in $\Sm$, i.e.\ $e^*f$ is a quotient of some
$\sheafify g$ in $\F/\sheafify X$.
\end{defn}

The following little observation will be used later on.

\begin{lem}\label{Slemepi}
The epis in $\F$ are precisely the sheafifications of epis in $\E$.
\end{lem}
\begin{prf}
First recall that epis in toposes are regular. Thus, since $\sheafify$ is a 
left adjoint it preserves regular epis. For the converse direction suppose
$e$ is an epi in $\F$. Consider its factorization $e = m \circ p$ in $\E$
where $m$ is monic and $p$ is an epi in $\E$. Then $e = \sheafify(m \circ p)
= \sheafify m \circ \sheafify p$ in $\F$. Since $\sheafify$ preserves monos 
and epis and $e$ is epic in $\F$ it follows that $\sheafify m$ is an iso. 
\end{prf}

Now we are ready to prove the main theorem of this subsection.

\begin{thm}\label{smallsub}
Suppose $\E$ is a topos with a natural numbers object $N$ and 
$\Sm$ is a class of small maps in $\E$ closed under power objects. 
If $\sheafify \dashv \incl : \F \inc \E$ is a subtopos then $\Sm_\F$ as
specified in Def.~\ref{SFdf} is a class of small maps in $\F$ 
which is closed under power objects.
\end{thm}

\begin{prf}
We will often (implicitly) use the fact that pullbacks in $\F$ preserve 
epis and maps in $\sheafify\Sm$.

This ensures for example that quasipullbacks of the form as considered in 
Def.~\ref{SFdf} are preserved by pullbacks along morphisms in $\F$.
Accordingly, it follows that $\Sm_\F$ is closed under pullbacks in $\F$,
i.e.\ validates condition {\bf (A0)}.

For showing that $\Sm_\F$ validates ${\bf (A1)}$ suppose that
\begin{diagram}[small]
B	\SEpbk	&\rCover		&  D       \\
\dTo^f	&			&  \dTo_g  \\
A		&\rCover_{p}		&  C       \\
\end{diagram} 
is a pullback in $\F$ where $f$ is in $\Sm_\F$ and $p$ is a cover in $\F$.
Since $f$ is in $\Sm_\F$ it fits into a quasipullback 
\begin{diagram}[small]
\sheafify Y & \rTo &  B  \\
\dTo^{\sheafify h} & \mathrm{qpb} & \dTo_f \\
\sheafify X & \rCover_e & A 
\end{diagram}
where $h$ is in $\Sm$ and $e$ is a cover in $\F$. 
Since quasipullbacks are closed under composition it follows that
\begin{diagram}[small]
\sheafify Y & \rTo &  B 	&  \rCover &  D\\
\dTo^{\sheafify h} &  & \dTo_f &  & \dTo_g  \\
\sheafify X & \rCover_e & A & \rCover_{p}		&  C
\end{diagram}
is a quasipullback. Thus, since $p \circ e$ is epic, it follows 
that $g$ is in $\Sm_\F$ as desired.

That $\Sm_\F$ validates condition {\bf (A2)} is immediate from the facts
that condition {\bf (A2)} holds for $\Sm$, that $\sheafify$ preserves $+$
and that $+$ preserves quasipullbacks.

That $\Sm_\F$ validates condition {\bf (A3)} is immediate from the fact 
that that $\sheafify$ preserves colimits and finite limits.

That $\Sm_\F$ validates {\bf (A4)}, i.e.\ that $\Sm_\F$ is closed under
composition, can be shown by adapting the proof of the analogous Lemma~2.15
of \cite{VDBM08}.

Obviously, $\Sm_\F$ validates condition {\bf (A5)} by its very definition
since quasipullbacks are closed under horizontal composition.

The proof that $\Sm_\F$ validates condition {\bf (A6)} is analogous to the
proof of case {\bf (A7)} of Proposition~2.14 of \cite{VDBM08}.

It is easy to check that {\bf (A7)} holds for $\Sm_\F$. Let $\pi$ be
a universal family for $\Sm$ then its sheafification $\sheafify\pi$ is 
universal for $\Sm_\F$ which can be seen by applying $\sheafify$ to
the respective diagram in the formulation of {\bf (A7)} and using the fact
that quasipullbacks are closed under horizontal composition.

Condition {\bf (A8)} holds for $\Sm_\F$ since sheafification preserves
natural numbers objects.

Condition {\bf (A9)} holds for $\Sm_\F$ since if $m$ is a mono in $\F$ 
then it is also a mono in $\E$ and thus by {\bf (A9)} for $\Sm$ we have
$m \cong \sheafify m$ is in $\Sm_\F$.

For showing that $\Sm_\F$ is closed under power objects one may adapt the
proof of Proposition~6.6 from \cite{VDBM08} proving an analogous result.
\end{prf}

\subsection{Small Maps in Subtoposes of $\RT(\A,\A_\#)$}

As a consequence of Theorem~\ref{smallsub} we obtain the following result.

\begin{thm}\label{subnestedthm}
Let $\Sm$ be the class of small maps in $\Asm(\A,\A_\#)$ as introduced 
in Theorem~\ref{Sclos1} and $\bar{\Sm}$ be the class of small maps in 
$\RT(\A,\A_\#)$ as introduced in Theorem~\ref{bennoextended}.
Suppose $\sheafify \dashv \incl : \E \inc \RT(\A,\A_\#)$ is a subtopos 
of $\RT(\A,\A_\#)$ induced by a closure operator $j$ on $\Trip(\A,\A_\#)$.
Then $\bar{\Sm}_\E$ as introduced in Theorem~\ref{smallsub} is a class of 
small maps in $\E$ closed under power objects and exponentiation.
\end{thm}
\begin{prf}
From Theorem~\ref{smallnesttop} we know that $\bar{\Sm}$ is a class of small
maps closed under power objects. Thus, we can apply Theorem~\ref{smallsub} from 
which it follows that $\bar{\Sm}_\E$ is a class of small maps in $\E$ which is 
closed under power objects and, accordingly, also under exponentiation.
\end{prf}

This result applies in particular to the subtoposes of $\RT(\A,\A_\#)$ as
considered in subsection~\ref{subtopex} and thus covers most of the examples
considered in van Oosten's book \cite{VOO08}.

\subsection{Models of $\IZF$ in Subtoposes of $\RT(\A,\A_\#)$}

From the main result of \cite{JM95} and our Theorem~\ref{subnestedthm} 
it follows that most of the toposes considered in \cite{VOO08} host models
of $\IZF$.

\begin{thm}\label{MainThm} There exist internal models for $\IZF$ in 
subtoposes of $\RT(\A,\A_\#)$ induced by local operators on $\Trip(\A,\A_\#)$.
\end{thm}

In case $\A = \A_\#$ due to \cite{KGVO05} we reobtain the realizability model 
for $\IZF$ as initially introduced by H.~Friedman in \cite{FRI73}, G.~Rosolini 
in \cite{ROS82} and D.~C.~McCarty in \cite{MCC84}. 

In case $\A = \A_\# = \Kcal_1$, the first Kleene algebra (corresponding to
number realizability), from Theorem~\ref{MainThm} it follows that the modified 
realizability topos $\Mod = \Mod(\Kcal_1) = \RT_m(\Kcal_1,\Kcal_1)$ 
from \cite{VOO97} hosts a model of $\IZF$. 
Thus, in $\IZF$ one cannot derive Markov's Principle from Church's Thesis.

\section{Conclusion}

Relying on the main result of \cite{JM95} we have shown that relative 
realizability toposes and modified relative realizability toposes host models 
of $\IZF$. In the unnested case, i.e.\ $\A = \A_\#$ we reobtain the well known 
realizability models for $\IZF$ and a modified realizability model for $\IZF$
which to our knowledge cannot be found in the existing literature. Moreover,
as pointed out to us by B.~van den Berg our results also show that his recent
Herbrand Realizability topos hosts a model of $\IZF$. 

We have obtained these new models for $\IZF$ in a quite uniform way 
using the methods of Algebraic Set Theory.
Of course, one could define in each single case these models of $\IZF$ in
a much more traditional and direct way. Using an appropriate adaptation of
the results in \cite{KGVO05} one can presumably show that these ``hand made'' 
models are equivalent to the ones we have obtained in this paper by more
abstract and general means.

\bibliographystyle{plain}

\end{document}